\documentclass{article}
\usepackage{amsfonts,amsmath,amssymb,amsthm,cite,enumerate,graphicx,latexsym,mathrsfs,color}
\numberwithin{equation}{section}

\newtheorem{definition}{\indent Definition}[section]

\newtheorem{theorem}{\indent Theorem}[section]

\title{{Oscillation of damped second order quasilinear wave equations with mixed arguments
\thanks{This work is supported in part by the National Natural Science Foundation of China (Grant No. 11671237)}}}

\author{{Ying Sui\thanks{Email: suiying4320@163.com}}~~~~{Huimin Yu\thanks{Corresponding author, Email: hmyu@amss.ac.cn}}\\
\small\textit{$$ School of Mathematics and Statistics, Shandong Normal University, Jinan 250014, China.}\\}

\begin{document}
\date{}
\maketitle

\begin{center}
\begin{minipage}{130mm}{\small
\textbf{Abstract}:
Following the previous work \cite{YingSui},  we investigate the impact of damping on the oscillation of smooth solutions to some kind of quasilinear wave equations with Robin and Dirichlet boundary condition. By using generalized Riccati transformation and technical inequality method, we give some sufficient conditions to guarantee the oscillation of all smooth solutions. From the results, we conclude that positive damping can ``hold back" oscillation.
At last, some examples are presented to confirm our main results.

}
\end{minipage}
\end{center}

\section{Introduction}

\ \ \ \ Following the work\cite{YingSui}, we continue to consider the oscillatory nature of quasilinear wave equation of the form
\begin{equation}\label{1.1}
\begin{aligned}
r(t)u^{\alpha-1}u_{tt}+p(x,t)u^{\alpha-2}u^2_t+ & \hat{p}(x,t)u^{\alpha-1}u_t+f(u,x,m(t))\\
&=a(t)\triangle u(x,t)+\sum_{k=1}^s a_k(t)\triangle u(x,\eta(t)),\ \ \
(x,t)\in G,
 \end{aligned}
\end{equation}
where $\alpha$ is the ratio of positive odd integers, $G=\Omega\times(t_0, \infty)$, $t_0>0$,
 $\Omega$ is a bounded domain in $R^N$ with a piecewise smooth boundary $\partial\Omega$,
 and $\triangle u(x,t)=\Sigma^N_{j=1}\frac{\partial^2u(x,t)}{\partial x_j^2}$.
 Here the first order term  $\hat{p}(x,t)u^{\alpha-1}u_t$ is added to consider the damping effect on the equation.
  Assumptions on the coefficients ($r, a, a_k$ and $p$) and functions ($f$, $ m$ and $\eta$) are the same as in \cite{YingSui}.
For convenient, we list them here:

 (H1)  $ r(t) \in C^1( (t_0,\infty),R_+)$; $  a(t), a_k(t)\in C( (t_0,\infty),R_+)$, $k\in I_s=\{1,2,\cdot\cdot\cdot, s\} $;
  $p(x,t)\in C(G,R_+)$  and $ p(x,t)\geq(\alpha-1)r(t)$ for $t>t_0$; $\hat{p}(x,t)\in C(G,R)$;

 (H2)   $ m(t)\in C((t_0,\infty),R_+)$, $ \eta(t)\in C((t_0,\infty),R_+)$, $m(t)$,  $\eta(t)$ is increasing,
  $m(t)\geq t$, and $\lim\limits_{t\rightarrow\infty}m(t)=\lim\limits_{t\rightarrow\infty}\eta(t)=\infty;$

 (H3)  $f\in C(R\times\Omega \times(t_0,\infty), R)$,
  $f(u,x,t) > q(x,t)u^\alpha(x,t)\geq q(t)u^\alpha(x,t)$ for some positive function $q(x,t)\in C(G,R_+)$ and $q(t)=\min_{x\in\bar{\Omega}}q(x,t)$.

 When $\hat{p}(x,t)\equiv 0$, the authors in\cite{YingSui} give three oscillation criteria for quasilinear wave equation (\ref{1.1})
  under the Robin boundary condition:
\begin{equation}\label{1.2}
\begin{aligned}
&\frac{\partial u(x,t)}{\partial \gamma}+\psi(x,t)u(x,t)=0,\ \ (x,t)\in \partial \Omega\times R_+,
 \end{aligned}
\end{equation}
and the Dirichlet boundary condition
\begin{equation}\label{1.3}
\begin{aligned}
&u(x,t)=0,\ \ (x,t)\in \partial \Omega\times R_+,
 \end{aligned}
\end{equation}
where $\gamma$ is the unit exterior normal vector to $\partial \Omega$ and $\psi(x,t)$ is a
 nonnegative continuous function on $\partial \Omega\times R_+$.  In this paper, we are interested in the influence of the additional (damped) term  $\hat{p}(x,t)u^{\alpha-1}u_t$ on the oscillation of equation (\ref{1.1})(\ref{1.2}) (or (\ref{1.1})(\ref{1.3})).

As we all know, oscillation phenomena is very common in our life and describes a form of material movement.
While damping is the characteristic that the oscillation amplitude decreases
gradually due to external effects or the inherent causes of the system itself, which plays a very important role in control theory, machinery, aviation, mechanics and other fields. In 2015, some oscillation theorems of certain nonlinear fractional partial differential equation with damping were established in \cite{3} by using differential inequality method as well as integral average method.
However, there are no deviating arguments considered.
 In 2016, the authors \cite{2} gave new characterizations on oscillation of third-order nonlinear damped delay differential equations.

In this paper, we consider the oscillation of the damped equation (\ref{1.1}) with advanced (or delayed) arguments.
To gain our aim, we first construct two new space variable integral average functions, and transform the quasilinear partial differential equation (\ref{1.1}) into an advanced functional differential equation with damping. Then some oscillation criteria are obtained for nonlinear problems (\ref{1.1})(\ref{1.2}) (or (\ref{1.1})(\ref{1.3})) by using the generalized Riccati transformation and technical inequality method.
From the results in this paper, we are convinced that positive damping can ``hold back" the  oscillation  in some sense.
 As far as we know, there is no any results on the oscillation of quasilinear damped wave equations with mixed arguments. We refer \cite{Li1,Li2,Yu,Cheng,HanZL,Ying,6,66,7,8} for more related results.\\

Before demonstrate our main results, we first give the definition of oscillation:
\begin{definition}\
 The solution $u(x,t)$ of problem (\ref{1.1}), (\ref{1.2}) (or (\ref{1.1}), (\ref{1.3}))  is said to be oscillatory in the domain $G$
 if for any positive number $\mu$, there exist $t>\mu$, $x_1 ,x_2\in \Omega$ such that $u(x_1 ,t) > 0$, $u(x_2 ,t) < 0$.
 Otherwise the solution $u(x ,t)$ is called nonoscillatory. If all of smooth solutions are oscillatory, we say the equation is oscillatory.
\end{definition}

\section{Main results}

\ \ \ \  Now, we present the main results of this paper.

\begin{theorem}\label{t5.2}
 Assume that  $(H1)$-$(H3)$ hold, $\hat{p}(x,t)\geq 0$ for any $(x,t)\in G$, $Q(t)=\frac{\alpha q(m(t))}{r(t)}$, $p_1(t)=\frac{\min_{x\in\Omega}\hat{p}(x,t)}{r(t)}$. If
 \begin{equation}\label{003.1}
\begin{aligned}
\limsup_{t\rightarrow\infty}\int_{t_0}^t e^{-\int^s p_1(\tau)d\tau}ds=\infty,
 \end{aligned}
\end{equation}
\begin{equation}\label{02.170}
\begin{aligned}
\limsup_{t\rightarrow\infty}\int_{t_0}^t sQ(s)ds=\infty,
 \end{aligned}
\end{equation}
and for sufficiently large $t\geq T\geq t_0$,
\begin{equation}\label{02.18}
\begin{aligned}
 m^{-1}(t)\int_{T}^{t}(s-T) Q(s)m(s)ds+(t-T)\int_{t}^\infty Q(s)ds>1.
 \end{aligned}
\end{equation}
Then  every smooth solution $u(x,t)$ of (\ref{1.1}) (\ref{1.2}) (or (\ref{1.1})(\ref{1.3}))  is oscillatory in $G$.
\end{theorem}

\textbf{Proof}.
By  (\ref{1.1}) and (H1) (H3), we have
\begin{equation}\label{02.101}
\begin{aligned}
&r(t)(u^{\alpha-1}u_{tt}+(\alpha-1)u^{\alpha-2}u^2_t)+ \hat{p}(x,t)u^{\alpha-1}u_t+q(m(t))u^\alpha(x,m(t))\\
<& a(t)\triangle u(x,t)+\sum_{k=1}^s a_k(t)\triangle u(x,\eta(t)).
 \end{aligned}
\end{equation}
Similar to the proof in \cite{YingSui}, we get
\begin{equation}\label{02.18}
\begin{aligned}
v''(t)+ p_1(t)v'(t)<-Q(t)v(m(t)),
 \end{aligned}
\end{equation}
where $v=\frac{1}{\alpha}\int_{\Omega}u^\alpha(x,t)dx$ for the Robin boundary condition, and
$v=\frac{1}{\alpha}\int_{\Omega}\varphi(x)u^\alpha(x,t)dx$ for the Dirichlet boundary condition.
Multiply both sides of (\ref{02.18}) by $e^{\int^tp_1(s)ds}$,
we get
$$
\left(v'(t)e^{\int^tp_1(s)ds}\right)'<-Q(t)v(m(t))e^{\int^tp_1(s)ds}<0,
$$
which implies that $v'(t)e^{\int^tp_1(s)ds}$ does not change sign eventually,
i.e.  $v'(t)$ does not change sign eventually,
 then there exists a $T\geq t_0$ such that either $v'(t) < 0$ or $v'(t) > 0$ for
any $t\geq T$.
If we assume that $v'(t)<0$, then
$$v(t)=v(T)+\int_{T}^tv'(s)ds\leq v(T)+v'(T)e^{\int^Tp_1(\tau)d\tau}\int_{T}^te^{-\int^sp_1(\tau)d\tau} ds,$$ for  any $t> T$,
 imply
$v(t)\rightarrow-\infty$ as $t\rightarrow\infty$, which is a
contradiction to  the positivity of  $v(t)$.
Therefore, $v'(t)>0$, for any $t\geq T$.
Then we have
\begin{equation}\label{02.180}
\begin{aligned}
v''(t)<- p_1(t)v'(t)-Q(t)v(m(t))\leq -Q(t)v(m(t))<0.
 \end{aligned}
\end{equation}
The following proof is similar to Theorem 2.1 in   \cite{YingSui}.
We can prove the theorem.\\

\begin{theorem}\label{t5.1}
Suppose $(H1)$-$(H3)$ hold, $Q(t)$, $p_1(t)$, and $\hat{p}(x,t)$ are the same as in Theorem 2.1.  Besides (\ref{003.1}) we
assume that there exist a constant $\beta$ and a function $\tau(t)\in C ((t_0,\infty),R_+)$ satisfying $\tau(t) \leq t$,
  $\tau'(t)\geq0$, $\tau''(t)\leq0$, such that the first order inequation
\begin{equation}\label{021}
\begin{aligned}
v'(t)-Q_1 (t)v(m(\tau(t)))>0
 \end{aligned}
\end{equation}
has no positive solution, where $ Q_1(t)=\int_t^{t+\beta} Q_0(s)ds$ and $ Q_0(t)=\min\{Q(t),Q(\tau(t))(\tau'(t))^2\}$. Then every smooth solution $u(x,t)$ of (\ref{1.1}) (\ref{1.2}) (or (\ref{1.1})(\ref{1.3})) is oscillatory in $G$.
\end{theorem}

\textbf{Proof}.
Base on the proof of Theorem 2.1, there exists $T>t_0$ such that $v'(t)>0$ and $v''(t)<-Q(t)v(m(t))$ for any $t\geq T$,
then the subsequent proof of this theorem is similar to Theorem 2.3 in \cite{YingSui}. We omit the details here.\\

\textbf{Remark 2.1}.\textit{ Compared the two theorems with Theorem 2.1 and Theorem 2.3 in \cite{YingSui}, we are convinced that
the additional assumption (\ref{003.1}) make oscillation more difficult. Specifically, when the damping coefficient $\hat{p}(x,t)\leq 0$, the assumption (\ref{003.1}) satisfied naturally, those two theorems are equivalent to the corresponding theorems in \cite{YingSui}. When $\hat{p}(x,t)>0$, the assumption (\ref{003.1}) really make things, which means positive damping can ``hold back" oscillation in some sense.}\\

 From the proof of Theorem 2.1, we notice that the positivity of $\hat{p}(x,t)$ make a very important role to get the inequality (\ref{02.180}).
Next, we give two other sufficient theorems in which the positivity condition of damping coefficient ($\hat{p}(x,t)>0$) is not necessary.
\begin{theorem}\label{t5.3}
 Suppose $(H1)$-$(H3)$ hold, $Q(t)$, $p_1(t)$, and $\hat{p}(x,t)$ are the same as in Theorem 2.1.  Besides (\ref{003.1}) we
assume that there exists a function $b(t)\in C^1( (t_0,\infty),R_+)$ such that
\begin{equation}\label{03.1}
\begin{aligned}
\limsup_{t\rightarrow\infty}\int_T^t \left(b(s)Q(s)-\left(\frac{b'(s)}{b(s)}- p_1(s)\right)^2\frac{b(s)}{4}\right)ds=\infty,\ \ \text{for\ \ some}\ \ T>t_0.
 \end{aligned}
\end{equation}
Then  every smooth solution $u(x,t)$ of (\ref{1.1}) (\ref{1.2}) (or (\ref{1.1})(\ref{1.3})) is oscillatory in $G$.
\end{theorem}

The proof of the above  theorem is similar to Theorem 2.2 in \cite{YingSui}. We omit the details here.\\

\begin{theorem}\label{t1}
 Assume that  $(H1)$-$(H3)$ hold, $Q_*(t)=\alpha q(m(t))$, $h(t)=\min_{x\in\Omega} \hat{p}(x,t)-r'(t)$.

 (1) If
 \begin{equation} \label{a}
 \limsup_{t\rightarrow\infty}\int_{t_0}^t  r^{-1}(s)e^{-\int^s\frac{h(\tau)}{r(\tau)}d\tau} ds=\infty,
 \end{equation}
and  there  exist  $T_1>t_0$,  $b\in C^1((t_0, \infty), R_+)$ such that
\begin{equation}\label{2.1}
\begin{aligned}
\limsup_{t\rightarrow\infty}\int^t_{T_1}\left[b(s)Q_*(s)-\frac{1}{4}r(s)b(s)\left(\frac{b'(s)}{b(s)}-\frac{h(s)}{r(s)}\right)^{2}\right]ds=\infty;
 \end{aligned}
\end{equation}

(2) If \begin{equation} \label{b}\limsup_{t\rightarrow\infty} \int_{t_0}^t  r^{-1}(s)e^{-\int^s\frac{h(\tau)}{r(\tau)}d\tau} ds<\infty,\end{equation}
  besides (\ref{2.1}) we further  assume that there exist $T_*\geq t_0$, $T_2\geq T_*$, and $T_3\geq T_*$ such that
\begin{equation}\label{2.1001}
\begin{aligned}
\limsup_{t\rightarrow\infty} \int_{T_3}^t\frac{ 1}{r(\tau)e^{\int^\tau\frac{h(s)}{r(s)}ds}} \int_{T_2}^\tau Q_*(s)\theta(m(s))e^{\int^s\frac{h(\eta)}{r(\eta)}d\eta}dsd\tau=\infty,
 \end{aligned}
\end{equation}
where $\theta(t)=\int_{t}^\infty\frac{1}{r(\tau)}e^{-\int^\tau_{T_*}\frac{h(s)}{r(s)}ds}d\tau$. Then every smooth solution $u(x,t)$ of (\ref{1.1})(\ref{1.2}) (or (\ref{1.1})(\ref{1.3})) is oscillatory in $G$.
\end{theorem}

\textbf{Proof}.  Suppose to the contrary that there is a nonoscillatory solution $u(x,t)$
of problem (\ref{1.1}) (\ref{1.2}) (or (\ref{1.1})(\ref{1.3})).
 Without loss of generality, we may assume that $u(x,t) \geq 0$, $u(x,m(t)) \geq 0$,
$u(x,\eta(t)) \geq0$ in $\Omega\times (t_1,\infty).$
Considering the left hand of (\ref{1.1}), we have
\begin{equation}\label{2.101}
\begin{aligned}
&\left[r(t)u^{\alpha-1}u_{tt}+r(t)(\alpha-1)u^{\alpha-2}u^2_t+r'(t)u^{\alpha-1}u_{t}\right]
+\left[\hat{p}(x,t)-r'(t)\right]u^{\alpha-1}u_{t}\\&+ q(m(t))u^\alpha(x,m(t))
< a(t)\triangle u(x,t)+\sum_{k=1}^s a_k(t)\triangle u(x,\eta(t)),
 \end{aligned}
\end{equation}
where we have used the facts that
$p(x,t)\geq (\alpha-1)r(t)$ in ($H_1$) and $f(u,x,m(t))> q(m(t))u^\alpha(x,m(t))$ in ($H_3$).

Similar the calculations in \cite{YingSui}, we set
$v_1=\frac{1}{\alpha}\int_{\Omega}u^\alpha(x,t)dx$ for the Robin boundary condition, and
$v_2=\frac{1}{\alpha}\int_{\Omega}\varphi(x)u^\alpha(x,t)dx$ for the Dirichlet boundary condition.
There yields $v_1$, $v_2$ satisfy the following inequality
\begin{equation}\label{2.7}
\begin{aligned}
\left(r(t)v'(t)\right)'+h(t)v'(t)+Q_*(t)v(m(t))<0,
 \end{aligned}
\end{equation}
for any $t>t_1>t_0$.

Next, we claim that there exists $t_2\geq t_1$ such that $v(t)$ is monotonic on $[t_2,+\infty)$.
In fact, from (\ref{2.7}), we have
$$
\left(r(t)v'(t)\right)'+h(t)v'(t)<-Q_*(t)v(m(t))\leq0.
$$
By setting $y(t) = r(t)v'(t)$, we get
\begin{equation}\label{2.7111}
\begin{aligned}
y'(t)+\frac{h(t)}{r(t)}y(t)<0.
 \end{aligned}
\end{equation}
Multiply both sides of (\ref{2.7111}) by $e^{\int^t\frac{h(s)}{r(s)}ds}$,
we get
\begin{equation}\label{2.190}
\begin{aligned}
\left(y(t)e^{\int^t\frac{h(s)}{r(s)}ds}\right)'<0,
 \end{aligned}
\end{equation}
which implies that $y(t)e^{\int^t\frac{h(s)}{r(s)}ds}$ is decreasing
and  $y(t)$ is eventually of one sign. Then $v'(t)$ has a fixed sign
for all sufficiently large $t$ and we have one of the following:
\begin{equation}\label{2.20}
\begin{cases}
\text{Case}\ (1): v'(t)>0;\\
\text{Case}\ (2): v'(t)<0.
 \end{cases}
\end{equation}

\textbf{Case (1)}.
 Define a generalized Riccati transformation by
$$
w(t)=\frac{b(t)r(t)v'(t)}{v(t)}>0.
$$
Since (\ref{2.7}) and $m(t)\geq t$, we get
\begin{equation}\label{2.22}
\begin{aligned}
w'(t)&=\left[\frac{b(t)r(t)v'(t)}{v(t)}\right]'
=r(t)v'(t)\left(\frac{b(t)}{v(t)}\right)'+(r(t)v'(t))'\frac{b(t)}{v(t)}\\
&< r(t)v'(t)\frac{b'(t)v(t)-b(t)v'(t)}{v^2(t)}-[h(t)v'(t)+Q_*(t)v(m(t))]\frac{b(t)}{v(t)}\\
&\leq\left(\frac{b'(t)}{b(t)}-\frac{h(t)}{r(t)}\right)w(t)-\frac{1}{r(t)b(t)}w^2(t)-b(t)Q_*(t)\\
&=-\left[\frac{1}{\sqrt{r(t)b(t)}}w(t)
-\frac{\sqrt{r(t)b(t)}}{2}\left(\frac{b'(t)}{b(t)}-\frac{h(t)}{r(t)}\right)\right]^2
+\frac{r(t)b(t)}{4}\left(\frac{b'(t)}{b(t)}-\frac{h(t)}{r(t)}\right)^2-b(t)Q_*(t)\\
&\leq\frac{1}{4}r(t)b(t)\left(\frac{b'(t)}{b(t)}-\frac{h(t)}{r(t)}\right)^{2}-b(t)Q_*(t).
 \end{aligned}
\end{equation}
Integration (\ref{2.22}) from $T_1>t_0$ to $t$, yields
$$
w(T_1)>\int^t_{T_1}\left[b(s)Q_*(s)-\frac{1}{4}r(s)b(s)\left(\frac{b'(s)}{b(s)}-\frac{h(s)}{r(s)}\right)^{2}\right]ds,
$$
which contradicts (\ref{2.1}).

\textbf{Case (2)}.
 By (\ref{2.190}), for $\tau\geq T_*\geq t_0$, we have
$$
y(T_*)e^{\int^{T_*}\frac{h(s)}{r(s)}ds}\geq y(\tau)e^{\int^\tau\frac{h(s)}{r(s)}ds},
$$
then
\begin{equation}\label{3.5}
\begin{aligned}
v'(\tau)\leq \frac{r(T_*)}{r(\tau)}v'(T_*)e^{-\int^\tau_{T_*}\frac{h(s)}{r(s)}ds}.
 \end{aligned}
\end{equation}
Integrating (\ref{3.5}) from  $t$ to  $\infty$,
\begin{equation}\label{3.6}
\begin{aligned}
 v(t)&\geq -r(T_*)v'(T_*)\int_{t}^\infty\frac{1}{r(\tau)}e^{-\int^\tau_{T_*}\frac{h(s)}{r(s)}ds}d\tau:=A\theta(t),
 \end{aligned}
\end{equation}
where $A=-r(T_*)v'(T_*)>0$, $\theta(t)=\int_{t}^\infty\frac{1}{r(\tau)}e^{-\int^\tau_{T_*}\frac{h(s)}{r(s)}ds}d\tau$.
Using (\ref{3.6}) in (\ref{2.7}), we find
\begin{equation}\label{3.7}
\begin{aligned}
-\left(r(t)v'(t)\right)'> h(t)v'(t)+Q_*(t)v(m(t))
\geq h(t)v'(t)+Q_*(t) A \theta(m(t)).
 \end{aligned}
\end{equation}
Define the function $V(t)=r(t)v'(t)<0$, and (\ref{3.7}) yields
\begin{equation}\label{3.8}
\begin{aligned}
V'(t)<- \frac{h(t)}{r(t)}V(t)- A Q_*(t)\theta(m(t)),
 \end{aligned}
\end{equation}
for any $t> T_*$. Multiply both sides of (\ref{3.8}) by $e^{\int^t\frac{h(s)}{r(s)}ds}$,
 we get
\begin{equation}\label{3.9}
\begin{aligned}
\left(V(t)e^{\int^t\frac{h(s)}{r(s)}ds}\right)'<-A Q_*(t) \theta(m(t))e^{\int^t\frac{h(s)}{r(s)}ds}.
 \end{aligned}
\end{equation}
Integrating (\ref{3.9}) from  $ T_2\geq T_*$ to  $t$,
\begin{equation}\label{3.10}
\begin{aligned}
V(t)e^{\int^t\frac{h(s)}{r(s)}ds}&< V(T_2)e^{\int^{T_2}\frac{h(s)}{r(s)}ds}- A \int_{T_2}^tQ_*(s)\theta(m(s))e^{\int^s\frac{h(\tau)}{r(\tau)}d\tau}ds\\
&\leq -  A\int_{T_2}^tQ_*(s) \theta(m(s))e^{\int^s\frac{h(\tau)}{r(\tau)}d\tau}ds,
 \end{aligned}
\end{equation}
then
\begin{equation}\label{3.11}
\begin{aligned}
v'(t)
< -\frac{ A }{r(t)e^{\int^t\frac{h(s)}{r(s)}ds}} \int_{T_2}^tQ_*(s) \theta(m(s))e^{\int^s\frac{h(\tau)}{r(\tau)}d\tau}ds.
 \end{aligned}
\end{equation}
Integrating (\ref{3.11}) from  $ T_3\geq T_*$ to  $\infty$, we get
$$
v(T_3)> A\int_{T_3}^\infty\frac{ 1}{r(\tau)e^{\int^\tau\frac{h(s)}{r(s)}ds}} \int_{T_2}^\tau Q_*(s) \theta(m(s))e^{\int^s\frac{h(\eta)}{r(\eta)}d\eta}dsd\tau,
$$
which is a contradiction to (\ref{2.1001}).\\

At last, following (\ref{2.20}), we claim that the assumption
 $ \int_{T}^t  r^{-1}(s)e^{-\int^s\frac{h(\tau)}{r(\tau)}d\tau} ds=\infty$ can odd out the Case (2) (i.e. $v'(t)<0$).
In fact, if $v'(t)<0$, then
$$v(t)=v(T)+\int_{T}^tv'(s)ds\leq v(T)+r(T)v'(T)e^{\int^T\frac{h(s)}{r(s)}ds}\int_{T}^t  r^{-1}(s)e^{-\int^s\frac{h(\tau)}{r(\tau)}d\tau} ds,\  \text{for\ any}\ t> T>t_0,$$
 imply
$v(t)\rightarrow-\infty$ as $t\rightarrow\infty$, which is a
contradiction to  the positivity of  $v(t)$.
Now, we finish the proof.

\section{Examples}

\textbf{Example 3.1}
  As an illustrative example, we consider the following equation
\begin{equation}\label{31}
\begin{aligned}
tu^{4}u_{tt}+u^{3}u^2_t
+u^{4}u_t+2u^5(x,2t)= \triangle u(x,t)+\sum_{k=1}^s (3+\cos kt)\triangle u(x,\frac{t}{2}),
 \end{aligned}
\end{equation}
with the  boundary conditions:
$\frac{\partial u(x,t)}{\partial \gamma}+tu(x,t)=0,\ \ (x,t)\in \partial \Omega\times R_+.$

Here $\alpha=5$,   $r(t)=t$, $p(x,t)=1$, $\hat{p}(x,t)=1$, $q(t)=1$, $m(t)=2t$,
$a(t)=1$, $a_k(t)=3+\cos kt$,  $\eta(t)=\frac{t}{2}$,
 $f(u,x,t)=2u^5(x,t)$, $b(t)=1$, then $Q_*(t)=5$.
It is easy to check that all hypotheses in case (1)
of Theorem 2.4 are satisfied, so we conclude the equation (\ref{31}) with
 Robin boundary condition is oscillatory.

\textbf{Example 3.2}
  We consider the following equation
\begin{equation}\label{32}
\begin{aligned}
t^2u^{2}u_{tt}+ uu^2_t
+2tu^{2}u_t+2t^4u^3(x,t+1)= \triangle u(x,t)+\sum_{k=1}^s (1+ kt)\triangle u(x,t+2),
 \end{aligned}
\end{equation}
with the  boundary conditions:
$u(x,t)=0,\ \ (x,t)\in \partial \Omega\times R_+.$

Here $\alpha=3$,   $r(t)=t^2$, $p(x,t)=1$, $\hat{p}(x,t)=2t$, $q(t)=t^4$, $m(t)=t+1$,
$a(t)=1$, $a_k(t)=1+kt$,  $\eta(t)=t+2$,
 $f(u,x,t)=2t^4u^3(x,t)$, $b(t)=1$.
It is easy to check that all hypotheses in case (2)
of Theorem 2.4 are satisfied, so the equation (\ref{32}) with Dirichlet boundary condition is oscillatory.

\end{document}